\documentclass[12pt]{article}
\usepackage{color}
\definecolor{darkblue}{rgb}{0.00,0.25,0.50}
\usepackage[colorlinks,filecolor=blue,citecolor=darkblue]{hyperref}

\usepackage{amsfonts,amssymb,amsmath,amsthm}
\usepackage{url}
\usepackage{enumerate}
\usepackage[ukrainian, russian, english]{babel}
\usepackage[cp1251]{inputenc}
\begin{document} \selectlanguage{ukrainian}
\thispagestyle{empty}

\title{}

UDC 517.51 \vskip 5mm

\begin{center}
\textbf{\Large Approximation of Nikolskii--Besov functional classes \boldmath{$S^{r}_{1,\theta}B(\mathbb{R}^d)$} by entire functions of a special form$^{\ast}$ }
\end{center}

\vskip 3mm

\begin{center}
\textbf{\Large  Наближення класів функцій Нікольського--Бєсова \boldmath{$S^{r}_{1,\theta}B(\mathbb{R}^d)$} цілими функціями спеціального вигляду }
\end{center}
\vskip0.5cm

\begin{center}
 S.~Ya.~Yanchenko\\ \emph{\small
Institute of Mathematics of NAS of
Ukraine, Kyiv}
\end{center}
\begin{center}
C.~Я.~Янченко \\
\emph{\small Інститут математики НАН України, Київ}
\end{center}
\vskip0.5cm

\begin{abstract}

We establish the exact-order estimates for the approximation of functions from the Nikol'skii--Besov classes $S^{\boldsymbol{r}}_{1,\theta}B(\mathbb{R}^d)$, $d\geqslant 1$, by entire function exponential type with some restrictions for their spectrum. The error of the approximation is estimated in the metric of the
Lebesgue space  $L_{\infty}(\mathbb{R}^d)$.

\vskip 5 mm

Одержано точні за порядком оцінки наближення функцій  з класів
 Нікольського--Бєсова $S^{\boldsymbol{r}}_{1,\theta}B(\mathbb{R}^d)$, $d\geqslant 1$, за допомогою цілих функцій експоненціального типу з певними обмеженнями
на їхній спектр.  Похибка наближення оцінюється у метриці простору Лебега $L_{\infty}(\mathbb{R}^d)$.

\end{abstract}

\vskip 0.7 cm

%%%%%% Вступ %%%%%%%%%%%%%%%%%%%%%%%%%%%%%%%%%%%%%%%%%%%%%%%%%%%%%%%%%%%%%%%%%%%%%%%%

\vskip 1mm  \textbf{1. Означення класів функцій та апроксимативних характеристик}.
У роботі продовжено (див. \cite{WangHeping_SunYongsheng_1995}--\cite{Yanchenko_YMG_2019}) вивчення апроксимативних характеристик класів функцій Нікольського--Бєсова
  $S^{\boldsymbol{r}}_{p,\theta}B(\mathbb{R}^d)$ у просторі $L_{\infty}(\mathbb{R}^d)$ в одновимірному та багатовимірному випадках і встановлено точні за
порядком оцінки наближення функцій із згаданих класів цілими функціями, зі
спектром зосередженим на множинах лебегова міра яких не перевищує $M$.

Простори функцій з домінуючою мішаною похідною $S^{\boldsymbol{r}}_{p,\theta}B(\mathbb{R}^d)$, які при значенні параметра $\theta=\infty$ тотожні просторам $S^{\boldsymbol{r}}_{p}H(\mathbb{R}^d)$ ($S^{\boldsymbol{r}}_{p,\infty}B(\mathbb{R}^d)\equiv S^{\boldsymbol{r}}_{p}H(\mathbb{R}^d)$), були введені
С.\,М.~Нікольським~\cite{Nikolsky_63} при $\theta=\infty$, а у випадку $1\leqslant \theta< \infty$~--- Т.\,І.~Амановим~\cite{Amanov_1965}. Означення даних  функціональних просторів у класичній формі були сформульовані С.\,М.~Нікольським та  Т.\,І.~Амановим  через кратні різниці та кратні модулі неперервності функцій.

\vskip 3mm

{\footnotesize \noindent $\overline{\ \ \ \ \ \ \ \ \ \ \ \ \ \ \ \ \ \ \ \ \ \ \ \ \ }$\\
\noindent  $^{\ast}$Публікація містить результати досліджень, проведених за грантом Президента України за конкурсним проектом №\,Ф84/177-2019 та виконана за часткової пiдтримки гранту НАН України дослідницьким лабораторiям/групам молодих вчених НАН України для проведення досліджень за пріоритетними напрямами розвитку науки i техніки, проект №04–02/2019.}

\newpage

Дамо означення просторів функцій Нікольського--Бєсова $S^{\boldsymbol{r}}_{p,\theta}B(\mathbb{R}^d)$ опосередковано через, так зване, декомпозиційне представлення елементів цих просторів. Уперше декомпозиційне представлення та відповідне йому нормування з'явилося у роботі С.\,М.~Нікольського та П.\,І.~Лізоркіна~\cite{Lizorkin_Nikolsky_1989} і, як з'ясувалося пізніше, відіграло ключову роль у дослідженнях, які пов'язані з апроксимацією класів функцій.  Це представлення істотно використовується для доведення одержаних результатів i базується на понятті перетворення Фур'є, яке можна означити використовуючи узагальнені функції  (див., наприклад, \cite{Lizorkin_69}], \cite[гл.~2]{Vladimirov}, \cite[гл.~11]{Berezanskyi2014}).

Наведемо спочатку необхідні означення та позначення, які у загальному будемо формулювати для $d\geqslant 1$ і в певних ситуаціях конкретизуємо їх для одновимірного випадку.

Нехай $\mathbb{R}^d$~--- $d$-вимірний евклідів простір з елементами
${\boldsymbol{x}=(x_1,\dots,x_d)}$ і
${(\boldsymbol{x},\boldsymbol{y})=x_1y_1+\dots+x_dy_d}$.  Нехай
${L_q(\mathbb{R}^d)}$, ${1\leqslant q\leqslant\infty}$,
--- простір вимірних на $\mathbb{R}^d$ функцій ${f(\boldsymbol{x})=f(x_1,\dots,x_d)}$
зі скінченною нормою
 $$
\|f\|_q:=
\Bigg(\int\limits_{\mathbb{R}^{d}}|f(\boldsymbol{x})|^{q}d\boldsymbol{x}
\Bigg) ^{\frac{1}{q}}, \ 1\leqslant q<\infty,
 $$
 $$
  \|f\|_{\infty}:=\mathop {\rm ess \sup}\limits_{\boldsymbol{x}\in \mathbb{R}^d}
  |f(\boldsymbol{x})|.
 $$

Далі нехай $S=S(\mathbb{R}^d)$~--- простір Л.~Шварца основних нескінченно
диференційовних на $\mathbb{R}^d$ комплекснозначних функцій
$\varphi$, які спадають на нескінченності разом зі своїми похідними
швидше за будь-який степінь функції $\left(x_1^2+\ldots+x_d^2\right)^{-\frac{1}{2}}$, що розглядається з відповідною топологією. Через $S'$ позначимо
простір лінійних неперервних функціоналів над $S$. Зазначимо, що
елементами простору $S'$ є узагальнені функції. Якщо $f\in S'$, $\varphi\in S$, то $\langle f,\varphi \rangle$  позначає значення $f$ на $\varphi$.

Через  $\mathfrak{F}\varphi$ та $\mathfrak{F}^{-1}\varphi$ будемо позначати відповідно пряме та обернене перетворення Фур'є функцій з просторів $S$ та $S'$.

Носієм неперервної на $\mathbb{R}^d$ функції $\varphi$ називається замикання множини точок ${\boldsymbol{x} \in \mathbb{R}^d}$, де $\varphi(\boldsymbol{x})\neq 0$, і позначається $\mbox{supp}\,\varphi$.

Узагальнена функція  $f$ перетворюється в нуль на відкритій множині $G$, якщо $\langle f,\varphi \rangle=0$ для всіх $\varphi \in S$ і $\mbox{supp}\,\varphi\subset G$. Об'єднання всіх околів, у яких $f$ перетворюється в нуль є відкритою множиною, яку називають нульовою множиною узагальненої функції $f$ і позначають $G_f$. Носієм узагальненої функції називають доповнення множини $G_f$ до $\mathbb{R}^d$, тобто замкнену множину $\mbox{supp}\,f=\bar{G}_f$.

Зазначимо, що кожна функція $f \in L_p(\mathbb{R}^d)$, $1\leqslant p \leqslant \infty$, визначає лінійний неперервний функціонал на $S$ згідно з формулою
$$
\langle f,\varphi \rangle = \int \limits_{\mathbb{R}^d} f(\boldsymbol{x})\varphi(\boldsymbol{x}) d\boldsymbol{x}, \ \ \varphi\in S,
$$
і, як наслідок, у цьому сенсі вона є елементом $S'$. Тому перетворення Фур'є функції ${f \in L_p(\mathbb{R}^d)}$, $1\leqslant p \leqslant \infty$, можна
розглядати, як перетворення Фур'є узагальненої функції $\langle f,\varphi \rangle$.

Далі нехай $K_m(t)=\int_{\mathbb{R}} k_m(\lambda)e^{-2\pi i \lambda t} d\lambda$, $m\in \mathbb{Z}_+$, $K_{-1}\equiv 0$,  де \\
\begin{minipage}{9 cm}
$$
 k_m(\lambda)=
 \begin{cases}
    1, & |\lambda|<2^{m-1}, \\
    2(1-\frac{|\lambda|}{2^m}), & 2^{m-1}\leqslant|\lambda|\leqslant2^{m}, \\
    0, &  |\lambda|>2^{m},
 \end{cases}
$$
\end{minipage}
\begin{minipage}{7 cm}
$$
 k_0(\lambda)=
 \begin{cases}
    1-|\lambda|, & 0\leqslant|\lambda|\leqslant1, \\
    0, & |\lambda|>1.
 \end{cases}
$$
\end{minipage}

Для кожного вектора ${\boldsymbol{s}=(s_{1},\dots,s_{d})}$, ${s_{j}\in
\mathbb{Z}_+}$, ${j=\overline{1,d}}$, покладемо
\begin{equation}\label{A_s_prod}
A^*_{\boldsymbol{s}}(\boldsymbol{x})=\prod\limits_{j=1}^d \big(K_{s_j}(x_j)-K_{s_j-1}(x_j)\big),
\end{equation}
$$
A^*_{\boldsymbol{s}}(f,\boldsymbol{x})=f(\boldsymbol{x})\ast A^*_{\boldsymbol{s}}(\boldsymbol{x})=
\int\limits_{\mathbb{R}^d}f(\boldsymbol{y}) A^*_{\boldsymbol{s}}(\boldsymbol{x}-\boldsymbol{y}) d\boldsymbol{y}.
$$
Також для $\boldsymbol{s}\in \mathbb{Z}^d_+$ розглянемо множину
$$
Q_{2^{\boldsymbol{s}}}^*=\big\{\boldsymbol{\lambda}=(\lambda_1,\dots,\lambda_d)\colon
\ \eta(s_j)2^{s_{j}-1}\leqslant |\lambda_j|<2^{s_j}, \lambda_j\in
\mathbb{R}, \ \ j=\overline{1,d}\big\},
$$
де $\eta(0)=0$ і $\eta(t)=1, \ t>0$ (відповідно $Q_{2^{s}}^*$ при $d=1$).

Справедливим є таке твердження.

\bf Лема А \rm  (див., наприклад, \cite{WangHeping_SunYongsheng_1999_AppT})\textbf{.} \it Нехай $1\leqslant p\leqslant\infty$, тоді для будь-якої функції ${f\in
L_p(\mathbb{R}^d)}$ маємо
 $$
f(\boldsymbol{x})=\sum\limits_{\boldsymbol{s}} A^*_{\boldsymbol{s}}(f,\boldsymbol{x})
 $$
 і, крім того, \rm $\mbox{supp} \, \mathfrak{F}A_{\boldsymbol{s}}(f,\boldsymbol{x})\subseteq Q_{2^{\boldsymbol{s}}}^*$. \rm

Зауважимо, що $A^*_{\boldsymbol{s}}(f,\boldsymbol{x})$ є аналогами ``блоків'' сум Валле Пуссена періодичних функцій багатьох змінних (див., наприклад, \cite{Temlyakov_1986m}).

 У прийнятих позначеннях простори $S_{p,\theta}^{\boldsymbol{r}} B(\mathbb{R}^d)$,
$1\leqslant p, \theta \leqslant \infty$, $\boldsymbol{r}>0$, можна означити таким
чином~(див., наприклад, \cite{WangHeping_SunYongsheng_1999_AppT}, \cite{WangHeping_1997_Q}):
$$
S^{\boldsymbol{r}}_{p,\theta}B(\mathbb{R}^d):=\Big\{f\in L_p(\mathbb{R}^d)\colon
\|f\|_{S^{\boldsymbol{r}}_{p,\theta}B}<\infty \Big\},
$$
де
\begin{equation}\label{Norm_dek_Sr1}
   \|f\|_{S^{\boldsymbol{r}}_{p,\theta}B(\mathbb{R}^d)}\asymp \Bigg(\sum \limits_{\boldsymbol{s}\geqslant
   0}2^{(\boldsymbol{s},\boldsymbol{r})\theta}
   \|A^*_{\boldsymbol{s}}(f,\cdot)\|_p^{\theta}\Bigg)^{\frac{1}{\theta}}
\end{equation}
  при $1\leqslant\theta<\infty$ і
\begin{equation}\label{Norm_dek_Sr1_inf}
   \|f\|_{S^{\boldsymbol{r}}_{p}H(\mathbb{R}^d)}\asymp \sup
   \limits_{\boldsymbol{s}\geqslant
   0} 2^{(\boldsymbol{s},\boldsymbol{r})}\|A^*_{\boldsymbol{s}}(f,\cdot)\|_p.
\end{equation}
\rm

Тут і надалі по тексту для додатних величин $A$ і  $B$ вживається запис  $A\asymp B$, який означає, що існують такі додатні  сталі $C_1$ та $C_2$, які не залежать від одного істотного параметра у величинах  $A$ і  $B$ (наприклад, у співвідношеннях (\ref{Norm_dek_Sr1}) і (\ref{Norm_dek_Sr1_inf})~--- від функції $f$), що ${C_1 A \leqslant B \leqslant C_2 A}$. Якщо тільки $B\leqslant C_2 A $ $\big(B \geqslant C_1 A\big)$, то пишемо
$B\ll A$ ${\big(B \gg A \big)}$. Всі сталі $C_i$, $i=1,2,\dots$, які
зустрічаються у роботі, залежать, можливо, лише від параметрів, що
входять в означення класу, метрики, в якій оцінюється похибка
наближення, та розмірності простору $\mathbb{R}^d$. Крім того нерівності типу $\boldsymbol{a}\leqslant \boldsymbol{b}$ $(\boldsymbol{a}>\boldsymbol{b})$ для векторів ${\boldsymbol{a}=(a_1,\dots,a_d)}$ та ${\boldsymbol{b}=(b_1,\dots,b_d)}$ будемо розуміти покоординатно, тобто $a_j\leqslant b_j$ $(a_j>b_j)$, $j=\overline{1,d}$. Також будемо використовувати записи $\boldsymbol{t}\geqslant 0$ $(\boldsymbol{t} > 0)$, якщо ${t_j\geqslant{0}}$ $(t_j> 0)$, $j=\overline{1,d}$, і $\boldsymbol{a}\neq \boldsymbol{b}$, якщо $a_i\neq b_i$ хоча б для одного $i$, $i=\overline{1,d}$.

Далі, замість
$S_{p,\theta}^{\boldsymbol{r}}B(\mathbb{R}^d)$ і
$S^{\boldsymbol{r}}_{p}H(\mathbb{R}^d)$  часто будемо використовувати  позначення
$S_{p,\theta}^{\boldsymbol{r}}B$ і $S^{\boldsymbol{r}}_{p}H$
відповідно ($S_{p,\theta}^{r}B$ і $S^{r}_{p}H$ при $d=1$).

У подальшому будемо вважати, що координати вектора
$\boldsymbol{r}=(r_1,\dots,r_d)$ впорядковані таким чином:
$0<r_1=r_2=\dots=r_{\nu}<r_{\nu+1}\leqslant\dots\leqslant r_d$.
Вектору $\boldsymbol{r}=(r_1,\dots,r_d)$ поставимо у відповідність вектор
$\boldsymbol{\gamma}=(\gamma_1,\dots,\gamma_d)$, $\gamma_j=r_j/r_1$, $j=\overline{1,d}$, а вектору $\boldsymbol{\gamma}$, в свою чергу,~--- вектор $\boldsymbol{\gamma}'$, де
$\gamma'_j=\gamma_j$, при $j=\overline{1,\nu}$ i
$1<\gamma_j'<\gamma_j$, при $j=\overline{\nu+1,d}$.

Окрім цього нагадаємо, що у випадку $1<p<\infty$ норму функцій з просторів $S_{p,\theta}^{\boldsymbol{r}} B(\mathbb{R}^d)$ можна означити в дещо іншій формі.

Нехай $A\subset \mathbb{R}^d$~--- деяка вимірна множина. Позначимо через
$\chi_{_A}$ характеристичну функцію множини $A$, і для ${f \in
L_p(\mathbb{R}^d)}$ покладемо
$$
\delta_{\boldsymbol{s}}^*(f,\boldsymbol{x})=\mathfrak{F}^{-1}(\chi_{Q_{2^{\boldsymbol{s}}}^*}\cdot
\mathfrak{F}f).
$$

Тоді простори $S_{p,\theta}^{\boldsymbol{r}} B$,
$1<p<\infty$, $1\leqslant  \theta \leqslant \infty$, $\boldsymbol{r}>0$, можна означити наступним
чином~\cite{Lizorkin_Nikolsky_1989}:
$$
S^{\boldsymbol{r}}_{p,\theta}B:=\Big\{f\in L_p(\mathbb{R}^d)\colon
\|f\|_{S^{\boldsymbol{r}}_{p,\theta}B}<\infty \Big\},
$$
де
\begin{equation}\label{Norm_dek_Sr}
   \|f\|_{S^{\boldsymbol{r}}_{p,\theta}B}\asymp \Bigg(\sum \limits_{\boldsymbol{s}\geqslant
   0}2^{(\boldsymbol{s},\boldsymbol{r})\theta}
   \|\delta_{\boldsymbol{s}}^*(f,\cdot)\|_p^{\theta}\Bigg)^{\frac{1}{\theta}}
\end{equation}
  при $1\leqslant\theta<\infty$ і
\begin{equation}\label{Norm_dek_Sr_inf}
   \|f\|_{S^{\boldsymbol{r}}_{p}H}\asymp \sup
   \limits_{\boldsymbol{s}\geqslant
   0} 2^{(\boldsymbol{s},\boldsymbol{r})}\|\delta_{\boldsymbol{s}}^*(f,\cdot)\|_p.
\end{equation}
\rm

Під класом $S^{\boldsymbol{r}}_{p,\theta}B$ будемо розуміти множину
функцій ${f \in L_p(\mathbb{R}^d)}$ для яких
$\|f\|_{S^{\boldsymbol{r}}_{p,\theta}B}\leqslant  1$, і при цьому збережемо
для класів $S^{\boldsymbol{r}}_{p,\theta}B$ ті ж самі позначення, що і для
просторів $S^{\boldsymbol{r}}_{p,\theta}B$.

Як видно з (\ref{Norm_dek_Sr1})\,--\,(\ref{Norm_dek_Sr_inf}), для $f\in S^{\boldsymbol{r}}_{p,\theta}B$, $1<p<\infty$, має місце співвідношення
\begin{equation}\label{As_deltas}
 \|\delta_{\boldsymbol{s}}^*(f,\cdot)\|_p\asymp \|A^*_{\boldsymbol{s}}(f,\cdot)\|_p.
\end{equation}

Перейдемо до означення апроксимативної характеристики, яка досліджується у роботі.

Нехай $\mathcal{L}\subset \mathbb{Z}^d_+$~--- деяка скінченна
множина,
${\mathfrak{M}=\mathfrak{M}(\mathcal{L})=\bigcup
\limits_{\boldsymbol{s}\in \mathcal{L}}Q^*_{2^{\boldsymbol{s}}}}$. Тоді
для ${f \in L_q(\mathbb{R}^d)}$, ${1\leqslant q \leqslant \infty}$, позначимо
$$
S_{\mathfrak{M}}(f,\boldsymbol{x})=\sum \limits_{\boldsymbol{s}\in
\mathcal{L}} \delta_{\boldsymbol{s}}^*(f,\boldsymbol{x}).
$$

Зауважимо, що $S_{\mathfrak{M}}(f,\boldsymbol{x})$ є цілою функцією,
яка належить простору $L_q(\mathbb{R}^d)$ і $\mbox{supp}\, S_{\mathfrak{M}}(f,x)\subseteq
\mathfrak{M}$.

Далі, для $f\in L_q(\mathbb{R}^d)$ розглянемо апроксимативну
характеристику
$$
 e_M^{\mathfrak{F}}\big(f\big)_q=\inf \limits_{\mathcal{L}\colon
 {\rm mes} \, \mathfrak{M} \leqslant M}
 \left\|f(\cdot)-S_{\mathfrak{M}}(f,\cdot)\right\|_q.
$$

Для $S^{\boldsymbol{r}}_{p,\theta}B(\mathbb{R}^d) \subset L_q(\mathbb{R}^d)$ позначимо
\begin{equation}\label{em-ort}
e_M^{\mathfrak{F}}\big(S^{\boldsymbol{r}}_{p,\theta}B\big)_q=\sup \limits_{f\in S^{\boldsymbol{r}}_{p,\theta}B} e_M^{\mathfrak{F}}\big(f\big)_q.
\end{equation}

Нагадаємо означення ще однієї апроксимативної характеристики, яку будемо використовувати.

Для $\boldsymbol{s}\in \mathbb{Z}^d_+$ означимо множину $Q_n^{\boldsymbol{\gamma}}$ таким чином:
$$
Q_n^{\boldsymbol{\gamma}}=\bigcup \limits_{(\boldsymbol{s},\boldsymbol{\gamma}) \leqslant n}Q_{2^{\boldsymbol{s}}}^*,
$$
де $n\in \mathbb{N}$. Множина $Q_n^{\boldsymbol{\gamma}}$ називається східчастим гіперболічним
хрестом і при цьому $\mbox{mes} \ Q_n^{\boldsymbol{\gamma}}\asymp
2^n n^{d-1}$ (див., наприклад, \cite{Lizorkin_Nikolsky_1989}), де
$\mbox{mes}\, Q_n^{\boldsymbol{\gamma}}$ позначає лебегову міру множини $Q_n^{\boldsymbol{\gamma}}$.

Крім того, для $f\in L_q(\mathbb{R}^d)$, $1\leqslant q \leqslant \infty$, покладемо
$$
S_{Q_n^{\boldsymbol{\gamma}}}(f,\boldsymbol{x})=\sum\limits_{(\boldsymbol{s},\boldsymbol{\gamma}) \leqslant n}\delta_{\boldsymbol{s}}^*(f,\boldsymbol{x}), \ \ \boldsymbol{x}\in  \mathbb{R}^d
$$
і означимо
\begin{equation}\label{EQN1}
\mathcal{E}_{Q_n^{\boldsymbol{\gamma}}}(f)_q=
\|f(\cdot)-S_{Q_n^{\boldsymbol{\gamma}}}(f,\cdot)\|_q \ \ \mbox{та} \ \ \mathcal{E}_{Q_n^{\boldsymbol{\gamma}}}(S^{\boldsymbol{r}}_{p,\theta}B)_q=\sup\limits_{f\in S^{\boldsymbol{r}}_{p,\theta}B}\mathcal{E}_{Q_n^{\boldsymbol{\gamma}}}(f)_q.
\end{equation}

Конкретизуємо означення величини $\mathcal{E}_{Q_n^{\boldsymbol{\gamma}}}(f)_q$ в
одновимірному випадку.

При $d=1$ кожна з множин $Q^*_{2^s}$  є об'єднанням напівінтервалів
${(-2^{s},-2^{s-1}]}$ та ${[2^{s-1},2^{s})}$, $s\in\mathbb{Z}_+$, з відповідною модифікацією при $s=0$. Тоді східчастий гіперболічний хрест, вироджується  в інтервал $(-2^n, 2^n)$, як об'єднання множин $Q^*_{2^s}$ для усіх $s\leqslant n$, $s\in\mathbb{Z}_+$, а саме $Q_n^1=\bigcup_{s\leqslant n}Q_{2^s}^*$. Крім того маємо $|Q_n^1| \asymp 2^n$, (де $|Q_n^1|$ позначає довжину інтервалу).

Означення величини (\ref{EQN1}) для $f\in L_q(\mathbb{R})$, $1\leqslant q \leqslant \infty$, можна переписати таким чином:
$$
\mathcal{E}_{2^n}(f)_q=\|f(\cdot)-S_{2^n}(f,\cdot)\|_q, \ \ \ \mathcal{E}_{2^n}(S^{{r}}_{p,\theta}B)_q=\sup\limits_{f\in S^{{r}}_{p,\theta}B}\mathcal{E}_{2^n}(f)_q,
$$
де
$$
S_{2^n}(f,x)=\sum \limits_{s \leqslant n}\delta^*_{2^s}(f,x)
$$
і
$$
\delta_{2^s}^*(f,x)=\mathfrak{F}^{-1}(
\chi_{Q^*_{2^s}}\cdot\mathfrak{F}f).
$$

Безпосередньо з означення апроксимативних характеристик (\ref{em-ort}) і (\ref{EQN1}) випливає, що у випадку
$\mbox{mes} \ Q_n^{\boldsymbol{\gamma}}\asymp \mbox{mes} \ \mathfrak{M}$ виконується співвідношення
\begin{equation}\label{eqn-em}
   e_M^{\mathfrak{F}}\big(S^{\boldsymbol{r}}_{p,\theta}B\big)_{q}\ll \mathcal{E}_{Q_n^{\boldsymbol{\gamma}}}\big(S^{\boldsymbol{r}}_{p,\theta}B\big)_{q}.
\end{equation}

Зауважимо, що величина (\ref{em-ort}) є неперіодичним аналогом найкращого ортогонального наближення, а величина (\ref{EQN1}) відповідає наближенню східчастою гіперболічною сумою Фур'є.  З дослідженнями  класів Нікольського--Бєсова періодичних функцій з домінуючою мішаною похідною, з точки зору знаходження порядкових оцінок різних апроксимативних характеристик, можна ознайомитися у монографіях В.\,М.~Темлякова~\cite{Temlyakov_1986m}, А.\,С.~Романюка~\cite{Romanyuk_2012m} і  D.~D\~{u}ng, V.\,N.~Temlyakov and T.\,Ullrich~\cite{Cross_2018}.

\textbf{2. Наближення функцій з класів \boldmath{$S^{r}_{1,\theta}B(\mathbb{R}^d)$} цілими функціями.} Перш ніж перейти до формулювання та доведення основних результатів, наведемо декілька допоміжних тверджень.

\bf Теорема А \rm \cite{Amanov_1965}\textbf{.} \it Нехай $1\leqslant  p, \theta \leqslant \infty$,
$1\leqslant  p\leqslant  q \leqslant \infty$ і  такий вектор $\boldsymbol{\rho}$, що
${\rho_j=r_j-\left(\frac{1}{p}-\frac{1}{q}\right)>0}$, $j=\overline{1,d}$. Тоді, якщо $f\in S_{p,\theta}^{\boldsymbol{r}}B(\mathbb{R}^d)$, то $f\in   S_{q,\theta}^{\boldsymbol{\rho}}B(\mathbb{R}^d)$
 і
 $$
  \|f\|_{S_{q,\theta}^{\boldsymbol{\rho}}B(\mathbb{R}^d)}\ll \|f\|_{S_{p,\theta}^{\boldsymbol{r}}B(\mathbb{R}^d)}.
 $$

\vskip 1 mm

\bf Теорема~Б \rm \cite[c.~150]{Nikolsky_1969_book}\textbf{.} \it
Якщо $1\leqslant p \leqslant q \leqslant \infty$, тоді для цілої
функції експоненціального типу $g_{\boldsymbol{\nu}}\in L_p(\mathbb{R}^d)$, $\boldsymbol{\nu}=(\nu_1,\ldots,\nu_2)$, $\nu_i\geqslant 0$, $i=\overline{1,d}$, має
місце нерівність різних метрик
$$
 \|g_{\boldsymbol{\nu}}\|_{L_{q}(\mathbb{R}^d)}\leqslant 2^d\left( \prod \limits_{j=1}^d
 \nu_k\right)^{\frac{1}{p}-\frac{1}{q}}\|g_{\boldsymbol{\nu}}\|_{L_p(\mathbb{R}^d)}.
$$\rm

\vskip 1 mm

 \bf Лема 1 \cite{Yanchenko_YMG_2019}. \it
Справедливою є така оцінка
 \begin{equation}\label{As_infty}
 \|A^*_{\boldsymbol{s}}(\cdot)\|_{\infty}\asymp 2^{\|\boldsymbol{s}\|_1},
 \end{equation}
 де $\|\boldsymbol{s}\|_1=s_1+\ldots + s_d$, $s_j\in \mathbb{Z}_+$, $j=\overline{1,d}$. \rm

\vskip 1 mm

 \bf Лема 2 \cite{Yanchenko_YMG_2019}. \it Нехай $1\leqslant p<\infty$, тоді має місце оцінка
 \begin{equation}\label{As_p}
 \|A^*_{\boldsymbol{s}}(\cdot)\|_{p}\asymp 2^{\|\boldsymbol{s}\|_1\left(1-\frac{1}{p}\right)}.
 \end{equation}\rm

\vskip 1 mm

\bf Лема 3 \cite{Yanchenko_YMG_2019}. \it
Справедливою є  така оцінка
 \begin{equation}\label{As_Sum_infty}
 \bigg\|\sum\limits_{(\boldsymbol{s},1)=n+1}A^*_{\boldsymbol{s}}(\cdot)\bigg\|_{\infty}\asymp 2^n n^{d-1}.
 \end{equation} \rm
\vskip 1 mm

Зауважимо, що оцінки лем~1--2, як слідує з їх доведення, виконуються і при $d=1$ з відповідною модифікацію. Наведемо також ще один результат, який будемо використовувати при доведенні одержаних результатів.

\vskip 3 mm  \bf Теорема~В \cite{Yanchenko_YMG_2019}. \it Нехай $r_1>1$, $1 \leqslant \theta \leqslant \infty$. Тоді виконується таке порядкове співвідношення
\begin{equation} \label{teor_1d}
  \mathcal{E}_{Q_n^{\gamma}}\big(S^{\boldsymbol{r}}_{1,\theta}B\big)_{\infty}\asymp
  2^{-n\left(r_1-1\right)}n^{(\nu-1)\left(1-\frac{1}{\theta}\right)}.
\end{equation}
\vskip 1 mm \rm

Справедливим є так твердження.

\vskip 1 mm \bf Теорема~1. \it Нехай $r>1$, $1 \leqslant \theta \leqslant \infty$. Тоді при $d=1$ виконується така порядкова оцінка
\begin{equation} \label{teor_em1infty_d1}
  e_M^{\mathfrak{F}}\big(S^{r}_{1,\theta}B(\mathbb{R})\big)_{\infty}\asymp  M^{-r+1}.
\end{equation}
\rm

\textbf{\textit{Доведення.}}  Перш за все зауважимо, що у випадку  $d=1$ оцінку (\ref{teor_1d}) можна записати таким чином:
\begin{equation} \label{teor_1d=1}
  \mathcal{E}_{2^n}\big(S^{r}_{1,\theta}B\big)_{\infty}\asymp
  2^{-n\left(r-1\right)}.
\end{equation}

Оскільки $r>1$, то згідно з теоремою~А існує таке число ${\rho}$, $\rho_j=r-1>0$, що для $f\in S_{1,\theta}^{{r}}B(\mathbb{R})$ маємо $f\in S_{\infty,\theta}^{{\rho}}B(\mathbb{R})$, тобто $f\in L_{\infty}(\mathbb{R})$.

Щоб встановити оцінки зверху в
(\ref{teor_em1infty_d1}) підберемо для $M$ число $n\in \mathbb{N}$
 із співвідношення $|Q_n^1|\leqslant M< |Q_{n+1}^1|$, тобто $M\asymp 2^n$. Тоді з
оцінки (\ref{teor_1d=1}), врахувавши (\ref{eqn-em}), будемо мати
$$
e_M^{\mathfrak{F}}\big(S^{r}_{1,\theta}B\big)_{\infty}\ll \mathcal{E}_{2_n}\big(S^{r}_{1,\theta}B\big)_{\infty}\asymp 2^{-n(r-1)}\asymp M^{-r+1}.
$$

Встановимо в (\ref{teor_em1infty_d1}) оцінку знизу.

Далі розглянемо функцію
$$
 f_1(x)=C_3 2^{-n r} A^*_{n}(x), \ \ C_3>0,
$$
тобто функція $f_1$ складається з одного ``блоку'' $A^*_{s}(x)$, який вибирається при $s=n$.

Покажемо, що $f_1$ належить класу
$S_{1,\theta}^{r}B$, $1 \leqslant \theta \leqslant \infty$. Оскільки, згідно з лемою~2 $\big\|A_{s}^*(\cdot)\big\|_1\asymp C_4$,  то
$$
 \|f_1\|_{S^{r}_{1,\theta}B}\asymp \left( \sum
 \limits_{s} 2^{s r \theta}
\|A_{s}^*(f_1,\cdot)\|_1^{\theta}\right)^{\frac{1}{\theta}}\asymp \left( 2^{n r \theta}  2^{-n r \theta} \right)^{\frac{1}{\theta}}=1
$$
при $1 \leqslant \theta < \infty$, і
$$
 \|f_1\|_{S^{r}_{1,\infty}}\asymp \sup
 \limits_{s}2^{sr}
\|A_{s}^*(f_1,\cdot)\|_1 \asymp 2^{n r} 2^{-nr} =1.
$$

Далі підібравши $M$  так, щоб  $|\widetilde{Q}_n^1|\leqslant 4M< |\widetilde{Q}_{n+1}^1|$, де $\widetilde{Q}_n^1=Q_{2^n}^*$, $|\widetilde{Q}_n^1| \asymp 2^n$,  і скориставшись лемою~1 можемо записати
$$
\|f_1(\cdot)-S_{\mathfrak{M}}(f_1,\cdot)\|_{\infty}\geqslant\big|\|f_1(\cdot)\|_{\infty}
-\|S_{\mathfrak{M}}(f_1,\cdot)\|_{\infty}\big|\gg
$$
$$
\gg 2^{-n r} (2^n -M)\gg
2^{-n r} \ 2^n \asymp M^{-r+1}.
$$

Оцінку знизу встановлено.

Теорему~1 доведено.\vskip 2 mm

\vskip 1 mm \bf Теорема~2. \it Нехай $r_1>1$, $1 \leqslant \theta \leqslant \infty$. Тоді при $d\geqslant2$ виконується така порядкова оцінка
\begin{equation} \label{teor_em1infty}
  e_M^{\mathfrak{F}}\big(S^{\boldsymbol{r}}_{1,\theta}B(\mathbb{R}^d)\big)_{\infty}\asymp
  \big(M^{-1} \log^{\nu-1} M \big)^{r_1-1} \big(\log^{\nu-1}M\big)^{1-\frac{1}{\theta}}.
\end{equation}
\rm

\vskip 1 mm

Перш ніж безпосередньо перейти до доведення теореми~2 зробимо зауваження. Оскільки $r_1>1$, то, згідно з теоремою~А, існує такий вектор $\boldsymbol{\rho}$, $\rho_j=r_j-1>0$, $j=\overline{1,d}$, що для $f\in S_{1,\theta}^{\boldsymbol{r}}B(\mathbb{R}^d)$ маємо $f\in   S_{\infty,\theta}^{\boldsymbol{\rho}}B(\mathbb{R}^d)$, тобто $f\in L_{\infty}(\mathbb{R}^d)$. Окрім того,
можемо стверджувати, що при деякому $1<q_0<\infty$, $f\in S^{\boldsymbol{\rho}}_{q_0,\theta}B$, де $\rho_j=r_j-\left(1-\frac{1}{q_0}\right)>0$, $j=\overline{1,d}$.

\vskip 1 mm

{\textbf{\textit{Доведення.}}} \ \rm  Оцінка зверху в
(\ref{teor_em1infty}), отримується з теореми~В. Оскільки ${\text{mes} \
Q_n^{\boldsymbol{\gamma}}\ll 2^n n^{\nu-1}}$, то підібравши для $M$
число $n\in \mathbb{N}$ із співвідношення ${\text{mes} \
Q_n^{\boldsymbol{\gamma}}\leqslant M< \text{mes} \
Q_{n+1}^{\boldsymbol{\gamma}}}$, тобто $M\asymp 2^n n^{\nu-1}$, з
оцінки (\ref{teor_1d}), будемо мати
$$
e_M^{\mathfrak{F}}\big(S^{\boldsymbol{r}}_{1,\theta}B\big)_{\infty}\ll
2^{-n\left(r_1-1\right)}n^{(\nu-1)\left(1-\frac{1}{\theta}\right)}\asymp
$$
$$
  \asymp\big(M^{-1}\log^{\nu-1}M\big)^{r_1-1}
  \big(\log^{\nu-1}M\big)^{\left(1-\frac{1}{\theta}\right)}.
$$

Перейдемо до встановлення оцінки знизу в (\ref{teor_em1infty}).
Попередньо зауважимо, що її достатньо отримати у випадку $\nu=d$.

Нехай
$$
 \Theta(n)=\big\{\boldsymbol{s}=(s_1,\ldots,s_d)\in \mathbb{Z}^d\colon \ s_1+\ldots+s_d=n \big \}
 \ \ \mbox{i} \ \ \widetilde{Q}_n=\bigcup
 \limits_{s\in\Theta(n) }Q_{2^{\boldsymbol{s}}}^*,
$$
тоді $\text{mes} \ \widetilde{Q}_n\asymp 2^n n^{\nu-1}$.

На відміну від одновимірного випадку, у залежності від значення параметра $\theta$ розглянемо функції
$$
 f_2(\boldsymbol{x})=C_5 2^{-n r_1} n^{-\frac{d-1}{\theta}} \sum
 \limits_{\boldsymbol{s}\in \Theta(n)} A^*_{\boldsymbol{s}}(\boldsymbol{x}), \ \ C_5>0,
$$
при $1\leqslant \theta<\infty$, i
$$
 f_3(\boldsymbol{x})=C_6 2^{-n r_1} \sum
 \limits_{\boldsymbol{s}\in \Theta(n)} A^*_{\boldsymbol{s}}(\boldsymbol{x}), \ \ C_6>0,
$$
якщо $\theta=\infty$.

Покажемо, що функції $f_2$ і $f_3$
належать класам
$S_{1,\theta}^{\boldsymbol{r}}B$ та $S_{1,\infty}^{\boldsymbol{r}}B$ відповідно. Оскільки, згідно з лемою~2, має місце оцінка $\big\|A_{\boldsymbol{s}}^*(\cdot)\big\|_1\asymp C_7$,  то
$$
 \|f_2\|_{S^{\boldsymbol{r}}_{1,\theta}B}\asymp \left( \sum
 \limits_{\boldsymbol{s}\in \Theta(n)}2^{(\boldsymbol{s},\boldsymbol{r})\theta}
\|A_{\boldsymbol{s}}^*(f_2,\cdot)\|_1^{\theta}\right)^{\frac{1}{\theta}}\asymp
$$
$$
\asymp 2^{-nr_1} n^{-\frac{d-1}{\theta}} \left( \sum
 \limits_{\boldsymbol{s}\in \Theta(n)}2^{(\boldsymbol{s},\boldsymbol{r})\theta}\| A_{\boldsymbol{s}}^*(\cdot)\|_1^{\theta}\right)^{\frac{1}{\theta}}\asymp
$$
$$
 \asymp 2^{-nr_1} n^{-\frac{d-1}{\theta}}\left( \sum
 \limits_{\boldsymbol{s}\in \Theta(n)}2^{r_1(\boldsymbol{s},1)\theta}\right)^{\frac{1}{\theta}}\ll
 n^{-\frac{d-1}{\theta}} \left( \sum
 \limits_{\boldsymbol{s}\in \Theta(n)}1\right)^{\frac{1}{\theta}}\ll 1.
$$

Для $f_3$ будемо мати
$$
 \|f_3\|_{S^{\boldsymbol{r}}_{1,\infty}}\asymp \sup
 \limits_{\boldsymbol{s}\in \Theta(n)}2^{(\boldsymbol{s},\boldsymbol{r})}
\|A_{\boldsymbol{s}}^*(f_3,\cdot)\|_1\asymp
$$
$$
 \asymp 2^{-nr_1} \sup \limits_{\boldsymbol{s}\in \Theta(n)}2^{(\boldsymbol{s},\boldsymbol{r})}
\|A_{\boldsymbol{s}}^*(\cdot)\|_1\asymp 2^{-nr_1} \sup \limits_{(\boldsymbol{s},1)=n+1}2^{(\boldsymbol{s},\boldsymbol{r})}
 \ll 1.
$$

Далі, через $\mathcal{L}'$ позначимо множину векторів  $\boldsymbol{s}$, таких що $\boldsymbol{s}\in \Theta(n)$, і щоб для множини ${\mathfrak{M}=\mathfrak{M}(\mathcal{L}')=\bigcup
\limits_{\boldsymbol{s}\in \mathcal{L}'}Q^*_{2^{\boldsymbol{s}}}}$ виконувалося співвідношення
\begin{equation}\label{Qn_4M}
\text{mes} \ \widetilde{Q}_n \leqslant 4M < \text{mes} \
\widetilde{Q}_{n+1},
\end{equation}
де $M=M(n)=\text{mes} \ \mathfrak{M}$.

Скориставшись лемами~1, 3 та співвідношення (\ref{Qn_4M}), врахувавши, що ${\text{mes} \ \widetilde{Q}_n\asymp 2^n n^{\nu-1}}$, можемо записати
$$
\|f_2(\cdot)-S_{\mathfrak{M}}(f_2,\cdot)\|_{\infty}\geqslant\big|\|f_2(\cdot)\|_{\infty}
-\|S_{\mathfrak{M}}(f_2,\cdot)\|_{\infty}\big|\gg
$$
$$
\gg 2^{-n r_1} n^{\frac{d-1}{\theta}}(2^n
n^{d-1}-M)\gg
2^{-n r_1} n^{\frac{d-1}{\theta}} 2^n
n^{d-1}=
$$
$$
=2^{-n(r_1-1)}n^{(d-1)\left(1-\frac{1}{\theta}\right)}
\asymp\big(M^{-1}\log^{d-1}M\big)^{r_1-1}
\big(\log^{d-1}M\big)^{\left(1-\frac{1}{\theta}\right)}.
$$

Аналогічно у випадку $\theta=\infty$,  отримаємо
$$
\|f_3(\cdot)-S_{\mathfrak{M}}(f_3,\cdot)\|_{\infty}\gg
\big(M^{-1}\log^{d-1}M\big)^{r_1-1}\log^{d-1}M.
$$

Оцінки знизу встановлено.

Теорему~2 доведено.\vskip 1mm

\vskip 2 mm

На завершення роботи зробимо деякі коментарі, щодо одержаних результатів.

Результати теореми~1 та 2 є новими і для для класів Нікольського  $S^{\boldsymbol{r}}_{1}H(\mathbb{R}^d)$, $d\geqslant 1$.

Як видно з одержаних теорем, оцінка величини $e_M^{\mathfrak{F}}\big(S^{r}_{1,\theta}B(\mathbb{R})\big)_{\infty}$ (теорема~1), не залежить від параметра $\theta$ на відміну від випадку  $d\geqslant2$ (теорема~2).

Порядкові оцінки величини $ e_M^{\mathfrak{F}}
\big(S^{\boldsymbol{r}}_{p,\theta}B\big)_q$ для  ряду інших
співвідношень між параметрами $p$, $q$ і $\theta$ встановлено в
\cite{Yanchenko_YMG_2010_8}, де також показано, що існують співвідношення між параметрами $p$, $q$, $\theta$
при яких величини $e_M^{\mathfrak{F}}
\big(S^{\boldsymbol{r}}_{p,\theta}B\big)_q$ і $
\mathcal{E}_{Q_n^{\gamma}}
\big(S^{\boldsymbol{r}}_{p,\theta}B\big)_q$ мають різні порядки.

На даний час також інтенсивно досліджуються й узагальнення класів Нікольського--Бєсова  з домінуючою мішаною похідною як періодичних так і неперіодичних функцій. У цьому напрямі відзначимо роботи М.\,М.~Пустовойтова~\cite{Pustovoitov_94}, \cite{Pustovoitov_99}, Sun Yongsheng, Wang Heping~\cite{SunYongsheng_WangHeping_1997}, С.\,А.~Стасюка~\cite{Stasjuk_04}, С.\,А.~Стасюка і С.\,Я.~Янченка~\cite{Stasuk_Yanchenko_Anal_math}, С.\,Я.~Янченка~\cite{Yanchenko_YMG_2010_1}, В.\,В.~Миронюка, С.\,Я.~Янченка~\cite{Yanchenko_Myronyuk_MS_2013}. \rm

Зауважимо, що у випадку $d=1$ класи Нікольського--Бєсова функцій з домінуючою мішаною похідною $S^{\boldsymbol{r}}_{p,\theta}B(\mathbb{R}^d)$ тотожні ізотропним та анізотропним класам Нікольського--Бєсова $B^{r}_{p,\theta}(\mathbb{R}^d)$ і $B^{\boldsymbol{r}}_{p,\theta}(\mathbb{R}^d)$. Знаходженню точних за порядком значень деяких апроксимативних  характеристик даних класів присвячені роботи  \cite{Yanchenko_YMG_2015}--\cite{Yanchenko_UMJ2018}.

\vskip 3 mm

\textbf{Contact information:}
Department of the Theory of Functions, Institute of Mathematics of National
Academy of Sciences of Ukraine, 3, Tereshenkivska st., 01024, Kyiv, Ukraine.

\vskip 3 mm

E-mail: \href{mailto:Yan.Sergiy@gmail.com}{Yan.Sergiy@gmail.com}

\end{document}